\documentclass[
a4paper,twoside,11pt,leqno]{amsart}

\usepackage[OT1]{fontenc}
\usepackage[utf8]{inputenc}
\usepackage{amsmath,amssymb,amsfonts,amsthm}
\usepackage[english]{babel}

\usepackage{hyperref}

\newcommand{\R}{\mathbb{R}}

\renewcommand{\d}{\,\mathrm{d}}

\newcommand{\abs}[1]{\left\lvert #1 \right\rvert}

\newcommand{\scalarprod}[1]{\left\langle #1\right\rangle}

\newcommand{\vect}[1]{\begin{pmatrix} #1 \end{pmatrix}}

\theoremstyle{plain}

\newtheorem*{prop*}{Proposition}

\newtheorem*{conj*}{Conjecture}
\theoremstyle{definition}

\theoremstyle{remark}

\numberwithin{equation}{section}

\title[The MOI for ISWs]{The Moment of Instability for\\Internal Solitary Waves}
\author{Andreas Klaiber}
\thanks{\noindent\emph{Present address:} Zentrum Mathematik, Technische Universit\"at M\"unchen, Boltzmannstr.\ 3, 85747
Garching, Germany. \emph{Phone:}  +49 89 289 17900. \emph{Email:} klaiber@ma.tum.de}
\keywords{Internal solitary waves, stratified fluids, moment of instability, Hamiltonian PDEs, incompressible Euler
equations}

\begin{document}
\maketitle
\begin{abstract}
  In this note, we define a moment of instability 
  $m(c)$ for internal solitary waves 
  in continuously stratified fluids,
  which seems not to have been done before.
  To underline the suitability of the proposed $m(c)$, we identify the relation $m''(c)=0$ as a formal Fredholm condition, and 
  we show that $m''(c)$ displays a definite sign for small-amplitude waves.
\end{abstract}
\section{Introduction}
Internal solitary waves (ISWs) are ecologically important since they are involved in mixing mechanisms and energy transport in
lakes and oceans \cite{Apel02,GOA04,PreussePeetersEtAl10,PreusseFreistuehlerEtAl12}. 
In this context, a widely used mathematical model consists of 
the 2D Euler equations for incompressible, inviscid fluids with non-constant density.
This model is given by the equations 
\begin{subequations}
\label{eq:euler-full}
\begin{align}
	\label{eq:euler-full-rho} \rho_t &= - u\rho_x - v\rho_y,\\
  \label{eq:euler-full-u} u_t &= - u u_x - v u_y -\frac{p_x}{\rho},\\
  \label{eq:euler-full-v} v_t &= - u v_x - v v_y -\frac{p_y}{\rho}-g,
\intertext{complemented by the incompressibility constraint}
  \label{eq:euler-incompr}   0 &= u_x + v_y,
\intertext{the boundary conditions}
	\label{eq:euler-full-boundarycondition}
	v(t,x,0) = 0 &\quad \text{and}\quad v(t,x,1) = 0,
\intertext{and the far-field conditions} 
(\rho,u,v,p)(t,\pm\infty,y)&=(\bar\rho(y),0,0,\bar p(y)),\quad 0\le y\le 1.
\end{align}
\end{subequations}
In \eqref{eq:euler-full}, density $\rho$, velocity  $(u,v)$, and pressure $p$ are functions of 
time $t$, horizontal position $x\in\R$ and vertical position $y\in[0,1]$, and  
the constant $g$ denotes acceleration due to gravity.
The far-field $(\bar\rho(y),0,0,\bar p(y))$, itself an $x$- and $t$-independent solution of   
\eqref{eq:euler-full-rho}-\eqref{eq:euler-full-boundarycondition} with 
\begin{multline*}
\bar\rho:[0,1]\to(0,\infty) \ \text{differentiable with}\ \ \bar{\rho}'(y)<0,\ 0\le y\le 1,\\
\text{and}\; \bar{p}(y)=-g\int_0^y \bar{\rho}(\eta)\d{\eta},
\end{multline*}
is called the \emph{quiescent state}. 
Travelling wave solutions 
\begin{equation*}
  (\rho,u,v,p)(t,x,y) = (\hat{\rho},\hat{u},\hat{v},\hat{p})(x-ct,y),\quad \text{with some $c>0$,}
\end{equation*}
of \eqref{eq:euler-full} are called \emph{internal solitary waves} (ISWs) of speed $c$; we
refer to \cite{Kirchgassner82,LankersFriesecke97,James97} for mathematical results on their existence.

In order to study the stability of ISWs, one could start from the linearization of \eqref{eq:euler-full} about a given ISW, as done
by the author 
in \cite{Klaiber13-Diss,Klaiber14} to find an Evans-function approach to stability. 

A different general approach to investigate the stability of solitary waves
is based on the moment of instability, see e.~g.
\cite{GrillakisShatahStrauss87}, and references therein.
Here, we want to establish the moment-of-instability (MOI) route to stability of ISWs.

\section{Definition of a moment of instability $m(c)$ for ISWs}
According to \cite{Benjamin84}, the Euler equations \eqref{eq:euler-full} for stratified fluids possess a Hamiltonian formulation.
In terms of the density $\rho$, the vorticity-like quantity $\sigma$, and the associated streamfunction $\psi$, defined
as the solution of  
\begin{equation*}
\sigma = - \nabla \cdot (\rho\,\nabla\psi), \quad \text{with}\quad \left.\psi\right|_{y=0,1} = 0,
\end{equation*}
this can be formulated as 
\begin{equation}
  \partial_t\vect{\rho\\ \sigma} = \mathcal{J}(\rho,\sigma) \, \left( \widetilde{\mathcal{H}}-c\widetilde{\mathcal{I}}
  \right)' (\rho,\sigma)
  \label{eq:euler-hamsyst}
\end{equation}
in the co-moving frame $t,\tilde{x}=x-ct,y$ with writing $x$ instead of $\tilde{x}$, 
where the Hamiltonian $\widetilde{\mathcal{H}}-c \widetilde{\mathcal{I}}$ is composed of the energy functional
\begin{equation}
  \widetilde{\mathcal{H}}(\rho,\sigma) = \int_\R\int_0^{1} \frac{1}{2}\rho \abs{\nabla\psi}^{2} + gy(\rho-\bar{\rho}) \d{y}\d{x}
  \label{eq:euler-H}
\end{equation}
and the momentum functional
\begin{equation}
  \widetilde{\mathcal{I}}(\rho,\sigma) = \int_\R \int_{0}^{1} y\sigma \d{y}\d{x},
  \label{eq:euler-I}
\end{equation}
and $\mathcal{J}=\mathcal{J}(\rho,\sigma)$ denotes the (state-dependent!) skew-symmetric operator 
\begin{equation}
  \mathcal{J}(\rho,\sigma) = 
  \begin{pmatrix} 
    0 & -\rho_x\partial_y+\rho_y\partial_x\\
		-\rho_x\partial_y+\rho_y\partial_x & -\sigma_x\partial_y+\sigma_y\partial_x
  \end{pmatrix}.
		\label{eq:euler-J}
\end{equation}
The Hamiltonian formulation \eqref{eq:euler-hamsyst}, however, does not directly yield a variational principle due to
the non-invertibility of $\mathcal{J}$. 
Concretely, as a stationary solution of \eqref{eq:euler-hamsyst} an ISW $(\rho^c,\sigma^c)$ satisfies
\begin{equation}
  0 = \mathcal{J}(\rho^c,\sigma^c) \, \left( \widetilde{\mathcal{H}}-c\widetilde{\mathcal{I}} \right)'(\rho^c,\sigma^c)
\end{equation}
but, as a little calculation reveals (see, e.~g., \cite[p.\ 35]{Benjamin84}),
\begin{equation}
  \left( \widetilde{\mathcal{H}}-c\widetilde{\mathcal{I}} \right)'\left(\rho^c,\sigma^c\right)
  = \vect{gy-\frac{1}{2}\abs{\nabla\psi^c}^2\\ \psi^c} \not= 0,
  \label{}
\end{equation}
i.~e., $(\rho^c,\sigma^c)$ is not a critical point of $\widetilde{\mathcal{H}}-c\widetilde{\mathcal{I}}$!

This issue can be overcome by modifying $\widetilde{\mathcal{H}}-c\widetilde{\mathcal{I}}$ without spoiling the Hamiltonian structure.
In fact, taking the quantities 
\begin{subequations}
  \label{eq:euler-casimir}
\begin{align}
  \Delta \mathcal{H}(\rho,\sigma) 
  &:= -\int_\R \int_{0}^{1} g\left\{\int_{\bar{\rho}(y)}^\rho  \bar{\rho}^{-1}(\varrho)\d{\varrho}\right\} \sigma
  \d{y}\d{x},\\
  \Delta \mathcal{I}(\rho,\sigma) 
  &:= -\int_\R \int_{0}^{1} \bar{\rho}^{-1}(\rho)  \sigma \d{y}\d{x},
\end{align}
\end{subequations}
it is easily verified that 
\begin{equation}
  \label{eq:condition-casimir}
  \mathcal{J}(\rho,\sigma)\left(\Delta\mathcal{H}-c\Delta \mathcal{I}\right)'(\rho,\sigma) = 0 \quad\text{and}\quad 
  \left(\mathcal{H} -c\mathcal{I}\right)'(\rho^c,\sigma^c) = 0
\end{equation}
with $\mathcal{H}:=\widetilde{\mathcal{H}}+\Delta\mathcal{H}$ and $\mathcal{I}=\widetilde{\mathcal{I}}+\Delta
\mathcal{I}$.
Therefore, replacing $\widetilde{\mathcal{H}}-c\widetilde{\mathcal{I}}$ with 
\begin{multline*}
\mathcal{H}-c\mathcal{I} \equiv 
  \int_\R\int_0^{1} \frac{1}{2}\rho \abs{\nabla\psi}^{2} + g
    \int_{\bar{\rho}(y)}^\rho \left\{ y-\bar{\rho}^{-1}(\varrho)\right\}\d{\varrho}\d{y}\d{x}\\
    \qquad
  -c\int_\R \int_{0}^{1} \left\{ y-\bar{\rho}^{-1}(\rho) \right\} \sigma \d{y}\d{x},
\end{multline*}
results in a modified Hamiltonian formulation such that ISWs are, indeed, critical points of the Hamiltonian.
This was already noticed by \cite{TurkingtonEydelandEtAl91,Abarbanel86} but, as far as the author is aware, has not
been used in connection with the stability of ISWs.
For background material on so-called Casimir functionals, for which $\Delta\mathcal{H}-c\Delta\mathcal{I}$ is an example, 
their systematic derivation and their use 
in hydrodynamic contexts, see \cite{Abarbanel86} and references therein.

%
Now, we are in a position to define the moment of instability for ISWs in the usual way:
\begin{equation}
  m(c) := \left( \mathcal{H}-c \mathcal{I} \right)(\rho^c,\sigma^c).
  \label{eq:moi-definition}
\end{equation}
Since $\left( \mathcal{H}-c \mathcal{I} \right)'(\rho^c,\sigma^c) = 0$ by construction, we immediately have the usual relation
\begin{equation*}
  m''(c) 
  \equiv \frac{\d^2{}}{\d{c}^2} \left( \mathcal{H}-c \mathcal{I} \right)(\rho^c,\sigma^c)
  = -\frac{\d{}}{\d{c}} \mathcal{I}(\rho^c,\sigma^c).
\end{equation*}

In the rest of the paper, we study this $m(c)$.
In Sec.\ 2 we show that $m''(c)<0$ for ISWs of sufficiently small amplitude.
In Sec.\ 3 we show in a quite general situation, which covers ours, that the condition $m''(c)=0$ can be read as a
formal Fredholm condition. 

\par\smallskip\textit{Acknowledgment.} 
I thank Heinrich Freist\"uhler for drawing my attention to the fact that
the meaning of $m''(c)=0$ can be understood without reference to a
possibly existing Evans function, and for stimulating discussions on this topic.

\section{Proving $m''(c)<0$ for small ISWs}
For small waves%
\footnote{We assume here the genericity condition $\int_0^1\bar{\rho}(y)\varphi_0^3(y)\d{y}\not=0$ which is necessary
for the validity of the approximate expressions; cf.\ \cite{KirchgassnerLankers93}.}%
, we have \cite{Benney66,James97}
\begin{align*}
  c &= c_0 + \varepsilon^2,\\
  \psi^c(x,y) &= \varepsilon^2 A(\varepsilon x)\varphi_0(y) + O(\varepsilon^4),\\
  \rho^c(x,y) &= \bar{\rho}(y) - \frac{1}{c_0} \varepsilon^2 A(\varepsilon x) \bar{\rho}'(y) \varphi_0(y) +
  O(\varepsilon^4),
\end{align*}
where 
\begin{equation*}
  A''(X) = -\frac{1}{s} A(X) - \frac{r}{s} A(X)^2
  \quad \text{and}\quad (\bar{\rho}(y)\varphi_0'(y))' = \frac{g}{c_0^2}\bar{\rho}'(y)\varphi_0(y).
\end{equation*}
With these expressions at hand, it is straightforward to evaluate $m''(c)$.

\begin{align*}
  \mathcal{I} (\rho^c, \sigma^c) &= \frac{1}{c} \int_\R\int_0^1 \rho^c\, \abs{\nabla\psi^c}^2 \d{y}\d{x}\\
  &= \frac{1}{c}\int_\R\int_0^1 
  \left(\bar{\rho}(y) - \frac{1}{c_0} \varepsilon^2 A(\varepsilon x) \bar{\rho}'(y) \varphi_0(y) +
  O(\varepsilon^4)\right)\\&\qquad \times \left( (\varepsilon^3A'(\varepsilon x)\varphi_0(y))^2+(\varepsilon^2 A(\varepsilon
  x)\varphi_0'(y))^2 + O(\varepsilon^5)\right) \d{y}\d{x}\\
  &= \frac{\varepsilon^4}{c_0} \int_\R\int_0^1 \bar{\rho}(y) A(\varepsilon x)^2 \varphi_0'(y)^2 
  \d{y}\d{x} + O(\varepsilon^5)\\
  &= \frac{\varepsilon^4}{c_0} \;\int_\R A(\varepsilon x)^2\d{x}\; \int_0^1 \bar{\rho}(y)\varphi_0'(y)^2\d{y}+ O(\varepsilon^5)\\
  &= \varepsilon^3\,\frac{1}{c_0}\; \int_\R A(X)^2\d{X} \;\int_0^1 \bar{\rho}(y)\varphi_0'(y)^2\d{y}+ O(\varepsilon^5)\\
  &= K\,(c-c_0)^{\frac{3}{2}} + O\left( (c-c_0)^{\frac{5}{2}} \right)
\end{align*}
with the finite, positive constant
\begin{equation*}
  K:=\frac{1}{c_0} \int_\R A(X)^2\d{X} \int_0^1 \bar{\rho}(y)\varphi_0'(y)^2\d{y}>0.
\end{equation*}
Hence, we derive that
\begin{equation*}
  m''(c) = -\frac{\d{}}{\d{c}} \mathcal{I}[(\rho^c,\sigma^c)] = 
  -\frac{3}{2}K \,(c-c_0)^{\frac{1}{2}} + O\left( (c-c_0)^{\frac{3}{2}} \right)  < 0
\end{equation*}
holds for $0\le c-c_0\ll 1$, i.~e., for sufficiently small waves.

\section{Characterizing $m''(c) = 0$ as a Fredholm condition}
To simplify the notation, we write $\phi=(\rho^c,\sigma^c)$ for the ISW in the following.
In the situation above, 
differentiating the profile equation
\begin{equation}
(\mathcal{H}-c\mathcal{I})'(\phi)=0 
\label{prof}
\end{equation}
with respect to the position yields 
\begin{equation}
  (\mathcal{H}-c\mathcal{I})''(\phi)\frac{\partial\phi}{\partial x}=0,
\label{dprofeqoverdx} 
\end{equation}
while differentiating it with respect to the speed results in
\begin{equation}
(\mathcal{H}-c\mathcal{I})''(\phi)\frac{\partial\phi}{\partial c}=\mathcal{I}'(\phi).
\label{dprofeqoverdc} 
\end{equation}
Eqs.\ \eqref{dprofeqoverdx}, \eqref{dprofeqoverdc} give 
\begin{equation}
\mathcal{J}(\mathcal{H}-c\mathcal{I})''(\phi)\frac{\partial\phi}{\partial x}=0
\label{trivialeigenfunction}
\end{equation}
and
\begin{equation}
\mathcal{J}(\mathcal{H}-c\mathcal{I})''(\phi)\frac{\partial\phi}{\partial c}=\mathcal{J}\mathcal{I}'(\phi).
\label{trivialhauptfunction}
\end{equation}
As
\begin{equation}
\mathcal{J}\mathcal{I}'(\phi)=-\frac{\partial\phi}{\partial x}, 
\label{JQT}
\end{equation}
eqs.\ \eqref{trivialeigenfunction} and \eqref{trivialhauptfunction}
state that $0$ is an at least double eigenvalue for 
$$\dot u= \mathcal{J}(\mathcal{H}-c\mathcal{I})'(u),$$
with $\frac{\partial\phi}{\partial x}$ 
as an eigenfunction and $\frac{\partial\phi}{\partial c}$ as 
a first-order generalized eigenfunction.\par\medskip
Now, a second-order generalized eigenfunction $\psi$ would solve 
\begin{equation}
\mathcal{J}(\mathcal{H}-c\mathcal{I})''(\phi)\psi=\frac{\partial\phi}{\partial c}.
\label{psi}
\end{equation}
According to the Fredholm alternative, eq.\ \eqref{psi} has a non-trivial
solution if and only if its right hand side $\frac{\partial\phi}{\partial c}$
is orthogonal to the solution $\chi$ of the adjoint homogeneous equation
\begin{equation}
0= (\mathcal{J}(\mathcal{H}-c\mathcal{I})''(\phi))^*\chi=-((\mathcal{H}-c\mathcal{I})''(\phi)\mathcal{J})\chi.
\end{equation}
As \eqref{dprofeqoverdx} and \eqref{JQT} imply 
\begin{equation}
0= -((\mathcal{H}-c\mathcal{I})''(\phi)\mathcal{J})\mathcal{I}'(\phi)\quad\text{and thus}\quad \chi=\mathcal{I}'(\phi),
\end{equation}
the existence of $\psi\neq 0$ consequently is equivalent to
\begin{equation}
0=\frac{d}{dc}\mathcal{I}(\phi)=\scalarprod{\mathcal{I}'(\phi),\frac{\partial\phi}{\partial c}} ,
\end{equation}
i.\ e., vanishing of the moment of instability.

\par\textit{Remarks.} (i) The above argument slightly varies the one given by Zumbrun in \cite{Zumbrun08}
(Sec.\ 1, between the statements of Corollary 1.4 and Remark 1.5).\\
(ii) This argument literally applies to the situation of Grillakis et.\ al.\ \cite{GrillakisShatahStrauss87} 
by changing to their notation 
\begin{equation*}
  E = \mathcal{H}, Q=\mathcal{I}, J=\mathcal{J}, \phi = (\rho^c,\sigma^c), T'(0)=\partial_x.
\end{equation*}
Hence, it can be applied to various contexts that fall into this class. 

\bibliographystyle{abbrv}
\bibliography{../../dfg/dfgbib}

%

\end{document}